\documentclass[12pt]{article}
\usepackage{amssymb}
\usepackage{amsmath}
\newtheorem{teo}{Theorem}

\newtheorem{defi}{Definition}

\newtheorem{rem}{Remark}
\newtheorem{lem}{Lemma}

\title{Uncertainty Principle for the Cantor Dyadic Group
}
\author{A. V. Krivoshein
\footnote{Faculty of Applied Mathematics and Control Processes, Saint Petersburg State University,
Universitetsky prospekt, 35, Peterhof, Saint Petersburg,   198504, 
Russia}, 
E. A. Lebedeva
\footnote{Mathematics and Mechanics Faculty, Saint Petersburg State University,
Universitetsky prospekt, 28, Peterhof,  Saint Petersburg,
 198504, Russia; St. Petersburg State Technical University, Department of Calculus, 
 Polytekhnicheskay 29, 195251, St. Petersburg, Russia}
}
\date{
krivosheinav@gmail.com, ealebedeva2004@gmail.com}

\begin{document}
\maketitle

\newcommand{\nul}{{\bf0}}
\newcommand{\rd}{{\mathbb R}^d}
\newcommand{\zd}{{\mathbb Z}^{d}}
\renewcommand{\r}{{\mathbb R}}
\newcommand{\z} {{\mathbb Z}}
\newcommand{\cn} {{\mathbb C}}
\newcommand{\n} {{\mathbb N}}

\begin{abstract}
We introduce a notion of  localization for dyadic functions, i.e. functions defined on the Cantor group. 
Localization is characterized by  functional $UC_d$ 
that is similar to the Heisenberg 
uncertainty constant for real-line 
functions. We are looking for  dyadic analogs of quantitative uncertainty principles. To justify our definition we use some test functions including dyadic scaling and wavelet functions. 
\end{abstract}



\section{Introduction}	

Good time-frequency localization of function $f: \mathbb{R}\to \mathbb{C}$ means that both function $f$ and its Fourier transform $Ff$ have sufficiently fast decay at infinity.  
The functional called the Heisenberg  uncertainty constant (UC) serves as a quantitative characteristic of this property.  Smaller UCs correspond to 
more localized functions. 
The uncertainty principle (UP) expresses a fundamental property of nature and can be stated as follows. If $f \neq 0$ then it is impossible for $f$ and $Ff$ to be 
sharply concentrated
simultaneously. 
In terms of the UC it means that there exists an absolute lower bound for the UC. 

There are numerous analogs and extensions  of this framework for different algebraic and topological structures.
For example, the localization of periodic functions is measured by means of the Breitenberger UC \cite{B}. For some particular cases of locally compact groups (namely a euclidean motion groups, non-compact semisimple Lie groups, Heisenberg groups) a counterpart of the UC is suggested in \cite{PrSi}. A generalization of operator interpretation for the UC is discussed in \cite{Selig}. These and many others related topics are described in the excellent survey \cite{FolSit}.   But to our knowledge, the question of a quantitative UC for the Cantor dyadic group has not been addressed in the literature.
In this paper we try to understand what "good localization" means for  functions defined on the  Cantor dyadic  group. So, a notion of the dyadic UC is suggested and justified. The existence of a lower bound is proven for the dyadic UC. We calculate this functional for dyadic scaling and wavelet functions and find good localized  dyadic wavelet frames.

We do not discuss qualitative UPs in this paper.   There exists a qualitative UP for a wide class of groups and the Cantor group belongs to the class 
(see p.224  (7.1) \cite{FolSit}). It is easy to see that  dyadic function $f_0=\chi_{[0,\,1)}=\widehat{f}_0,$ where $\widehat{f}$ is a Walsh-Fourier transform of $f$ (see the definition in Section \ref{Not}), satisfies the extremal equality in this UP.     There are a lot of results in this direction (see \cite{HJ}, \cite{Gro} and references therein).

The paper is organized as follows. First, we introduce  necessary notations and auxiliary results.
In section \ref{locsec}, we formulate  the  definition of the dyadic UC, prove a dyadic UP,  answer the question how to calculate the dyadic UC    in some particular important cases. 
In section \ref{examsec}, we calculate the dyadic UC for Lang's wavelet and looking for wavelet frames with  small dyadic UCs.

\section{Notations and Auxiliary Results}
\label{Not}
Let 
$
x=\sum_{j \in \mathbb{Z}}x_j 2^{-j-1}
$
 be a dyadic expansion of $x \in [0,\,\infty)=\mathbb{R}_+,$ where $x_j \in \{0,\,1\}.$ 
For $x=p 2^n,$ $p \in \mathbb{N},$ $n \in \mathbb{Z},$ there are two possible expansions, one terminates in 0's and another does in 1's. We choose the first one, that is  $x_j \to 0$ as $j \to \infty.$ 
\texttt{The dyadic sum} of $x$  and $y$ is defined by 
$$
x \oplus y :=\sum_{j \in \mathbb{Z}}|x_j-y_j| 2^{-j-1}.
$$
 Then $[0,\,\infty)$ is metrizable with the distance between $x,\,y$ defined to be $x \oplus y$.
 A function that is continuous from the $\oplus$-topology to the usual topology is called
 \texttt{W-continuous}.
It is well known (see \cite[sections 1.3, 9.1]{SWS}, \cite[sections 1.1, 1.2]{GES}) that this framework  is a representation of 
\texttt{the Cantor dyadic group}, i.e. the Cartesian product of countably many copies of $\mathbb{Z}_2$, the discrete cyclic group of order 2 (the set $\{0,\,1\}$ with discrete topology and modulo 2 addition).   

\texttt{The Walsh-Fourier transform }of $f\in L_1(\mathbb{R}_+)$ is defined by 
$$
\widehat{f}(t):= \int_{\mathbb{R}_+} f(x) {\rm w}(t,\,x)\,dx, 
$$
 where
 the function ${\rm w}(t,\,x):=(-1)^{\sum_{j \in \mathbb{Z}}t_j x_{-j-1}}$ is the representation for a character of the dyadic group. 
 The Walsh-Fourier transform inherits many properties from the Fourier transform  (see \cite[sections 9.2, 9.3]{SWS}). For example,  the Plancherel theorem holds 
 $$
 \int_{\mathbb{R}_+}f(x)\,\overline{g(x)}\,dx \,=\, 
 \int_{\mathbb{R}_+}\widehat{f}(x)\,\overline{\widehat{g}(x)}\,dx,
 $$
 for   $f,\,g,\,\widehat{f},\,\widehat{g} \in L_1(\mathbb{R}_+)$ with standard extension to $L_2(\mathbb{R}_+).$ 
 Functions ${\rm w}(n,x),$ where $n = 0,\,1,\,2, \dots$ are called \texttt{the Walsh functions}.   They form an orthonormal basis for $L_2([0,\,1))$. The Walsh system is a dyadic analog of the trigonometric system. 
  

 \texttt{The fast Walsh-Fourier transform }of $x=(x_k)_{k=\overline{0,2^n-1}}\in \r^{2^n}$ is defined by 
$
	c=x{\rm W},
$
 where 
 ${\rm W}=2^{-\frac n2}({\rm w}(m,\,k/2^n))_{k,m=0}^{2^n-1}=\{\bar{\omega}^{n}_{k,m}\}_{k,m=0}^{2^n-1}$
 is \texttt{the normalized Walsh matrix} 
(see \cite[section 9.7]{SWS} accurate within the normalization).
 The matrix ${\rm W}$ is orthogonal, symmetric, and unitary ${\rm W}^{-1}={\rm W}.$
 
The concept of a dyadic derivative is quite different from its classical counterpart (see \cite[section 1.7]{SWS}, \cite[section 6.3]{SA}). 
The function 
$$
f^{[1]}(x):=\sum_{j\in\mathbb{Z}}2^{j-1}(f(x)-f(x\oplus 2^{-j-1}))
$$
is called \texttt{the dyadic derivative} of $f$ at $x$.
The inherited properties are the following
$$
{\rm w}^{[1]}(n,\,x)\,=\, n {\rm w}(n,\,x), \quad \quad
\widehat{f^{[1]}}(t)\,=\,t \widehat{f}(t).
$$
But unfortunately the dyadic derivative does not support some natural properties such as the chain rule and the rule $(f g)'\,=\,fg'+f'g.$

Let $H$ be a separable Hilbert space. If there exist  constants $A,\,B>0$
such that for any $f \in H$ the following inequality holds
$
A \|f\|^2 \leq \sum_{n=1}^{\infty} \left|(f,\,f_n)\right|^2 \leq B \|f\|^2,
$ 
then the sequence $(f_n)_{n \in \mathbb{N}}$ is called  \texttt{a frame}  for $H.$
If $A=B \, (=1),$ then the sequence $(f_n)_{n \in \mathbb{N}}$ is called \texttt{a (normalized) tight frame} for $H.$

If the set of functions $\psi_{j,k}(x):=2^{j/2}\psi(2^j x  \oplus k)$ forms a frame or a basis of $L_2(\mathbb{R}_+),$ then it is called \texttt{a dyadic wavelet frame or basis}. Using the routine procedure, it can be generated from multiresolution analysis starting with an auxiliary function, that is a scaling function $\varphi$.


The foundation of the dyadic  (Walsh) analysis is contained in \cite{SWS}, \cite{GES}. 
The concept of a dyadic wavelet function and elements of  multiresolution analysis  theory for the Cantor dyadic group is developed in \cite{Lang} and later in \cite{FP},
\cite{F12}.

\section{Localization of Dyadic Functions}
\label{locsec}

The quantitative characteristic of the time-frequency localization is the uncertainty constant (UC). 
Originally, the concept of an uncertainty constant and principle was introduced for the real line case in 1927. 
\texttt{The Heisenberg uncertainty constant}  of $f \in L_2(\mathbb{R})$ is the functional 
$UC_H(f):=\Delta_{f}\Delta_{Ff}$ such that
\begin{equation*}
\label{intUC}
\begin{array}{ll}
\Delta_{f}^2:=\frac{1}{\|f\|^{2}_{L^2(\mathbb{R})}}\int_{\mathbb{R}}(x-x_{f})^2|f(x)|^2\,d x, &
\Delta_{Ff}^2:=\frac{1}{\|Ff\|^{2}_{L^2(\mathbb{R})}}
 \int_{\mathbb{R}}(t-t_{Ff})^2|Ff(t)|^2\,d t, \\
\end{array}
\end{equation*}
$$
\begin{array}{ll}
x_{f}:=\frac{1}{\|f\|^{2}_{L^2(\mathbb{R})}}\int_{\mathbb{R}}x|f(x)|^2\,d x, &
t_{Ff}:=\frac{1}{\|Ff\|^{2}_{L^2(\mathbb{R})}}\int_{\mathbb{R}}
t|Ff(t)|^2\,d t, \\
\end{array}
$$  
where $Ff$ denotes the Fourier transform of $f.$
It is well known that 
$UC_H(f)\geq 1/2$ for a function $f \in L_2(\mathbb{R})$ and the minimum is attained on the Gaussian.  Let us make some preliminary remarks to motivate the definition of a localization characteristic for the dyadic case. 

\begin{rem}
\label{rem1}
It is easy to see that $x_f$ is the solution of the minimization problem
$$
\min_{\tilde{x}}\int_{\mathbb{R}}(x-\tilde{x})^2|f(x)|^2\,d x.
$$  
Hence, the squared $UC_H$ takes the form
$$
\frac{1}{\|f\|^{2}_{L^2(\mathbb{R})}} \min_{\tilde{x}}\int_{\mathbb{R}}(x-\tilde{x})^2|f(x)|^2\,d x\ \ 
\frac{1}{\|Ff\|^{2}_{L^2(\mathbb{R})}} \min_{\tilde{t}}\int_{\mathbb{R}}(t-\tilde{t})^2|Ff(t)|^2\,d t.
$$
\end{rem}

\begin{rem}
\label{rem2}
It is well known that $x_f$ equals to the integral mean value of the function $f,$ while 
$\Delta_{f}$ means the dispersion with respect to the    $x_f$. The sense of the sign "-" in the definition of $\Delta_{f}$ is the distance between $x$ and $x_f.$ Thus, we have
$$
UC^2_H=\frac{1}{\|f\|^{2}_{L^2(\mathbb{R})}}
\min_{\tilde{x}}\int_{\mathbb{R}}{\rm dist}^2(x,\,\tilde{x})|f(x)|^2\,d x
\ \ 
\frac{1}{\|Ff\|^{2}_{L^2(\mathbb{R})}} \min_{\tilde{t}}\int_{\mathbb{R}}{\rm dist}^2(t,\,\tilde{t})|Ff(t)|^2\,d t.
$$   
\end{rem}

Now we are ready to introduce the definition of a localization characteristic for the dyadic setup.

\begin{defi}
\label{def1}
Suppose $f\in L_2(\mathbb{R}_+)$ is a complex valued dyadic function, then the functional 
$$
UC_d(f):=V(f)V(\widehat{f}), \quad \mbox{ where }
$$
$$
V(f):=
\frac{1}{\|f\|^2_{L_2(\mathbb{R}_+)}} 
\min_{\tilde{x}}\int_{\mathbb{R}_+}(x \oplus \tilde{x})^2|f(x)|^2\,d x,
$$
$$
V(\widehat{f}):=
\frac{1}{\|\widehat{f}\|^2_{L_2(\mathbb{R}_+)}}
\min_{\tilde{t}}\int_{\mathbb{R_+}}(t \oplus \tilde{t})^2|\widehat{f}(t)|^2\,d t
$$
is called \texttt{the dyadic uncertainty constant (the dyadic UC)}  of the function $f.$
\end{defi}

\begin{rem}
\label{minimum}
Suppose  $g$ is a bounded dyadic complex-valued function, $g(x),\,x g(x) \in L_2(\mathbb{R}_+)$. 
 We denote $G(y):=\int_{\mathbb{R}_+}(x \oplus y)^2|g(x)|^2\,d x$. 
 Since  $g(x),\,x g(x) \in L_2(\mathbb{R}_+)$ and $x\oplus y < x+y$ it follows that $G(y)$ is finite for 
 $y\in \mathbb{R}_+$
 Then there exists a point $y^{\ast}$ such that
   $\min_{y} G(y)= G(y^{\ast}).$
   Indeed, it is clear that $y^{\ast}$ can not be outside the interval $[0,\,2^n)$ for some probably large $n\in\mathbb{N}$ depending on $g$. It can be checked that  $[0,\,2^n)$ is compact in the dyadic topology. The function $x\oplus y$ is $W$-continuous, therefore $G$     is $W$-continuous.  It is well known that under these conditions, the image $G([0,\,2^n))$ is compact. Finally, since  $G([0,\,2^n))\subset \mathbb{C},$ it follows that 
   $G([0,\,2^n))$ is bounded and closed.
\end{rem}

\textbf{Example 1.}
\label{ex1}
 Let $\chi_M$ be a characteristic function of a set $M.$
Denote $f_1(x)= \chi_{[0,\,1/4)}(x)$ and $g_1(x)= \chi_{[3/4,\,1)}(x).$ Then it is easy to calculate their Walsh-Fourier transforms 
$\widehat{f_1}=\chi_{[0,\,4)}/4$
and
$\widehat{g_1}={\rm w}\left(3,\,\cdot/4\right)\chi_{[0,\,4)}/4.$    
It is natural to characterize "the dispersion" of these functions  by means of 
the diameters of their supports. Thus, ${\rm diam}[0,\,1/4):=\sup_{x,\,y\in[0,\,1/4)}(x\oplus y)=1/4,$
 ${\rm diam}[3/4,\,1)=1/4,$ and  ${\rm diam}[0,\,4)=4.$  So, these functions should have the same localization. On the other side, let us consider the functions $f_2(x)=\chi_{[0,\,3/8)}(x)$ and $g_2(x)= \chi_{[3/4,\,9/8)}(x).$ Their Walsh-Fourier transforms are 
 $\widehat{f_2}=\chi_{[0,\,4)}/4+{\rm w}\left(1,\,\cdot/4\right)\chi_{[0,\,8)}/8$ 
 and
 $\widehat{g_2}={\rm w}\left(3,\,\cdot/4\right)\chi_{[0,\,4)}/4+{\rm w}(1,\,\cdot)\chi_{[0,\,8)}/8.$
 Calculating the  diameters we get  
${\rm diam}[0,\,3/8)=1/2,$  ${\rm diam}[3/4,\,9/8)=2,$ and 
${\rm diam}[0,\,8)=8.$
 So, the first function should be more localized.
Indeed, Table \ref{tab1} shows that our suppositions are correct. 
Columns named $\tilde{x}_0(f)$ and $\tilde{t}_0(f)$ mean  sets of $\tilde{x}$ and $\tilde{t}$ minimizing the functionals $\int_{\mathbb{R}_+}(x \oplus \tilde{x})^2|f(x)|^2\,d x$ and $\int_{\mathbb{R}_+}(t \oplus \tilde{t})^2|\widehat{f}(t)|^2\,d t$ respectively.

\begin{table}[ht]
\centering
\caption{The dyadic uncertainty constants: Example 1.}
\begin{tabular}{ccccccc} 
\hline
 $f$ &  $\|f\|^2(=\|\widehat{f}\|^2)$
  & $\tilde{x}_0(f)$
  & $\tilde{t}_0(f)$
  & $V(f)$
  & $V(\widehat{f})$
  & $UC_d(f)$   \\
   \hline
$f_1$ & $1/4$ & $[0,\,1/4)$ & $[0,\,4)$   &  $1/48$ & $16/3$    & $1/9$ \\
$g_1$ & $1/4$ &  $[3/4,\,1)$ & $[0,\,4)$   &  $1/48$ & $16/3$  & $1/9$ \\
$f_2$ & $3/8$ & $[0,\,1/8)$ & $[0,\,2)$   & $3/64$ & $8$  & $3/8$ \\
$g_2$ & $3/8$ &  $[3/4,\,7/8)$ & $[0,\,4)$   &  $71/64$ & $32/3$  & $71/6$\\
 \hline
\end{tabular}
\label{tab1}
\end{table}

\begin{rem}
\label{rem3}
The operator interpretation of the UC does not work for the  
dyadic setup.
Let $P$ and $M$ be self-adjoint, symmetric or normal operators defined on a Hilbert space,
$[P,\,M]_{-}:=PM-MP$ be a \texttt{commutator} of $P$ and $M$, and $[P,\,M]_{+}:=PM+MP$ be an \texttt{anticommutator} of $P$ and $M.$ The following inequality named \texttt{the Schr\"oedinger uncertainty principle} (see \cite{Sch}) is a simple consequence of the Cauchy-Bunyakovski-Schwarz inequality
\begin{equation*}
\label{SchrodUP}
	\|Mf-\beta f\|^2  \|Pf-\alpha f\|^2\geq 
	\frac{1}{4}\left(|([P,M]_{-}f,\,f)|^2+\left|([P,M]_{+}f,\,f)-2\alpha \beta  \|f\|^2\right|^2\right),
\end{equation*}
where
$
\beta:=(Mf,\,f)/\|f\|^2, \ \alpha:=(Pf,\,f)/\|f\|^2.
$
 It gives two functionals both used as the UCs: the first one is more traditional, but some authors (see \cite{Selig})  exploit the second one as well
 \begin{equation}
UC_{-}(f):=\frac{\|Mf-\beta f\|  \|Pf-\alpha f\|}{ |([P,M]_{-}f,\,f)|} \geq 1/2 
\label{UC-}
\end{equation}
\begin{equation}
UC_{+}(f):=\frac{\|Mf-\beta f\|  \|Pf-\alpha f\|}{ \left|([P,M]_{+}f,\,f)-2\alpha \beta  \|f\|^2\right|} \geq 1/2.
\label{UC+}
\end{equation}
   Defying in (\ref{UC-}) $Pf(x)=i\,f'(x)$ and $Mf(x)=x\,f(x),$ one get the Heisenberg UC in $L_2({\mathbb R})$.
   The dyadic extension of this framework has the following trouble.  If the  inner product $(P_Hf,\,M_Hf)$ is real-valued then the mean value of the commutator $([P,\,M]_{-}f,\,f)=2 i\, \Im(P_Hf,\,M_Hf)$
   vanishes.
  In classical setup the inner product  is pure imaginary for a real-valued $f.$ 
   But for natural choice of dyadic operators on $L_2({\mathbb R_+})$, namely $Pf(x)=f^{[1]}(x)$ and $Mf(x)=x f(x)$, it turns out to be real-valued. Thus, one get identical zero in the denominator of (\ref{UC-}). The reason of the trouble is the difference between the operators $i\, f'$ and $ f^{[1]}.$ It is caused by the definitions of respective characters and the properties of derivatives, namely $({\rm e}^{i\,t})'=i\,{\rm e}^{i\,t}$ and $({\rm w}(n,\,t))^{[1]}=n {\rm w}(n,\,t),$ the imaginary unit appears only in the classical case.      
  
  A dyadic counterpart of (\ref{UC+}) does not give an adequate characteristic of localization. Indeed, it equals to infinity for the very well localized function $f_0:=\chi_{[0,\,1)},$ $\widehat{f}_0=f_0,$ while,
    $UC_d(f_3)=1/9.$ 
\end{rem}

There is a lower bound for $UC_d,$ so we get an uncertainty principle for the dyadic Cantor group.
\begin{teo}
\label{UP}
For any function 
$f\in L_2(\mathbb{R}_+),$ the following inequality holds
$$
UC_d(f)\geq C,
\mbox{ where }
C \simeq 8.5 \times 10^{-5}.
$$
\end{teo}  

\textbf{Proof.}
Suppose $f_1(x):={\rm w}(\tilde{t},\,x) f(x\oplus \tilde{x}),$ then 
$\widehat{f}_1(t):={\rm w}(t,\, \tilde{x}) \widehat{f}(t\oplus \tilde{t})$ and
it is straightforward calculation to see that 
$$
\int_{\mathbb{R_+}}(t \oplus \tilde{t})^2|\widehat{f}(t)|^2\,d t=
\int_{\mathbb{R_+}} t^2|\widehat{f_1}(t)|^2\,d t,
$$
\begin{equation}
\label{simp}
\int_{\mathbb{R_+}}(x \oplus \tilde{x})^2|f(x)|^2\,d x=
\int_{\mathbb{R_+}} x^2|f_1(x)|^2\,d x.
\end{equation}
So, it is sufficient to prove 
$$
\|x g(x) \|  \,  \|t \widehat{g}(t)\|\geq \sqrt{C} \|g\|^2.
$$
It can be done in the same manner as its classical counterpart (see \cite[Theorem 1.1, Corollaries 1.2, 1.3]{PrSi}).
\begin{enumerate}
	\item Let $E$ be a measurable subset of $\mathbb{R}_+,$ $|E|$ be a Lebesgue measure of $E$, and $0<\theta<1/2.$ 
Then 
$$
\left(\int_E |\widehat{f}|^2\right)^{1/2} \leq K_1(\theta) |E|^{\theta} \|x^{\theta}f(x)\|_2,
\mbox{ where  }
K_1(\theta)=(2 \theta)^{-2 \theta}(1-2\theta)^{\theta-1}.
$$
Indeed, suppose $B=[0,\,b),$  $B'=[b,\,\infty).$
Then 
$
\left(\int_E |\widehat{f }|^2\right)^{1/2} \leq 
\left(\int_E |\widehat{f \chi_B}|^2\right)^{1/2} +
\left(\int_E |\widehat{f\chi_{B'}}|^2\right)^{1/2}.
$
Using definition of the Walsh-Fourier transform, the Cauchy-Bunyakovskii-Schwarz inequality,  and elementary properties of integrals we get for the first and the second summands
\begin{gather*}
\left(\int_E |\widehat{f \chi_B}|^2\right)^{1/2} \leq
|E|^{1/2} \sup_{E} |\widehat{f\chi_{B}}| \leq
|E|^{1/2} \|f \chi_B\|_1 \leq 
|E|^{1/2} \|x^{-\theta}\chi_B(x)\|_2 \|x^{\theta} f(x)\|_2 
\\
=
|E|^{1/2} (1-2\theta)^{-1/2} b^{-\theta+1/2} \|x^{\theta} f(x)\|_2,
\\
\left(\int_E |\widehat{f\chi_{B'}}|^2\right)^{1/2} \leq
\|f \chi_B\|_2 \leq
\sup_{B'} x^{-\theta} \|x^{\theta} f(x)\|_2 \leq
b^{-\theta} \|x^{\theta} f(x)\|_2.
\end{gather*}	 
So,
	\begin{gather*}
	\left(\int_E |\widehat{f}|^2\right)^{1/2} \leq
	\left(|E|^{1/2} (1-2\theta)^{-1/2} b^{-\theta+1/2}+b^{-\theta}\right)\|x^{\theta} f(x)\|_2.
	\end{gather*} 
	It remains to minimize the right side over $b$ 
	($b_{\min} = 4\theta^2 |E|^{-1} (1-2\theta)^{-1}$) 
	to get the desired inequality.  
	\item 
	Let us prove 
	$
	\|f\|_2^2 \leq 2 K_1(\theta)\|x^{\theta} f(x)\|_2 \|t^{\theta} \widehat{f}(t)\|_2 
	$
	for $0<\theta<1/2.$
	Denote $E=[0,\,r),$ $E'=[r, \infty).$
	Then using the first item, we obtain
	\begin{gather*}
	\|f\|^2_2=\|\widehat{f}\|^2_2=
	\int_E |\widehat{f}|^2 + \int_{E'} |\widehat{f}|^2 \leq
	K_1^2(\theta) r^{2\theta} \|x^{\theta} f(x)\|_2^2+ r^{-2 \theta}\|t^{\theta} \widehat{f}(t)\|_2^2.
	\end{gather*} 
Minimizing the last expression over $r$ ($r_{\min}=\|t^{\theta} \widehat{f}(t)\|_2^{1/(4\theta)}(K_1^2(\theta)\|x^{\theta} f(x)\|_2)^{-1/(4\theta)}$)
we get the necessary inequality.	
	
	\item Since the function $g(\alpha):=\left(\|x^{\alpha} f(x)\|_2 \|f\|^{-1}_2\right)^{1/\alpha}$ decreases for $\alpha >0$ ($g'_{\alpha} > 0$), then 
	$$
	\|x^{\alpha} f(x)\|_2 \leq \|f\|_2^{1-\alpha/ \beta} \|x^{\beta} f(x)\|^{\alpha/\beta}_2  
	$$ 
	for $0<\alpha <\beta.$
	
	\item Applying the last inequality ($\alpha = \theta$) to  item 2 we obtain 
	$$
	\|f\|^2_2\leq 2 K_1(\theta) \|x^{\theta} f(x)\|_2 \|t^{\theta} \widehat{f}(t)\|_2 \leq
	2 K_1(\theta) \|f\|_2^{2-2\theta/ \beta}  \|x^{\beta} f(x)\|_2^{\theta/\beta} \|t^{\beta} \widehat{f}(t)\|_2^{\theta/\beta},
	$$
	thus
	$$
	\|f\|^2_2\leq (2 K_1(\theta))^{\beta/\theta} \|x^{\beta} f(x)\|_2 \|t^{\beta} \widehat{f}(t)\|_2. 
	$$
	So, choosing $\beta = 1$ we have 
	$$
	\|x f(x)\|_2 \|t \widehat{f}(t)\|_2 \geq C(\theta) \|f\|^2_2,
\mbox{ where }
C(\theta)=(2 K_1(\theta))^{-1/\theta}.
	$$
	To get the dyadic uncertainty principle it remains to maximize $C^2(\theta)$ over $\theta$, $\max_{\theta}C^2(\theta) \simeq C^2(0.382)\simeq 8.5\times  10^{-5}.$
	$\Box$
\end{enumerate}

It is not easy to calculate $UC_d$ for an arbitrary function because of the dyadic minimization problem underlying in the definition of $UC_d.$    The following result gives a possible way to calculate the dyadic UC on a wide class of functions. The minimization problem adds up to exhaustive search among $2^n$ variants.

\begin{lem}
\label{th1}
Let $f(x)=\chi_{[0,\,1)}(x) \sum_{k=0}^{\infty}a_k {\rm w}(k, x)$
be a uniformly convergent series restricted on $[0,\,1),$ 
 $f_n(x)=\chi_{[0,\,1)}(x) \sum_{k=0}^{2^n-1}a_k {\rm w}(k, x)$
be its partial sum, $V(f)<+\infty, $ $V(\widehat{f})<+\infty.$ Then the dyadic UC takes the form
$$
UC_d(f)\,=\,\lim_{n \to \infty}V(f_n)V(\widehat{f}_n), \mbox{ where }
$$
$$
V(f_n)=\frac{\min_{k_0=\overline{0,2^n-1}}\sum_{k=0}^{2^n-1} |c_{k\oplus k_0}|^2 ((k+1)^3-k^3)2^{-2n}/3}{\sum_{k=0}^{2^n-1}|a_k|^2},
$$
\begin{equation*}
\label{finiteV}
V(\widehat{f}_n)=\frac{\min_{k_1=\overline{0,2^n-1}}\sum_{k=0}^{2^n-1} |a_{k\oplus k_1}|^2 ((k+1)^3-k^3)/3}{\sum_{k=0}^{2^n-1}|c_{k}|^2},
\end{equation*}
and
$
c:=(c_{k})_{k=\overline{0,2^n-1}} 
$
is the fast Walsh-Fourier transform of $a:=(a_{k})_{k=\overline{0,2^n-1}}$. 
\end{lem}

\textbf{Proof.}
Suppose $\Delta_{k,n}:=[k 2^{-n},\,(k+1)2^{-n}),$
$ k=0,\dots,2^n-1,$  $n= 0,1,\dots$ is a dyadic interval,
 $\xi_{k,n}:=\chi_{\Delta_{k,n}}$ 
is the characteristic function of $\Delta_{k,n},$ and 
	$f_n(x)=\sum_{k=0}^{2^n-1} b_k \xi_{k,n}(x)$ is a representation of $f_n$ with respect to the orthogonal system 
	$\{\xi_{k,n},:\  k=0,\dots,2^n-1, \, n= 0,1,\dots\}.$ It is easy to find a connection between 
	$a=(a_k)_{k=\overline{0,2^n-1}}$ and $b=(b_k)_{k=\overline{0,2^n-1}}.$ Indeed,
	$$
	\sum_{k=0}^{2^n-1}a_k {\rm w}(k, x)=f_n(x)=\sum_{k=0}^{2^n-1}b_k \xi_{k,n}(x).
	$$
	The Walsh-Fourier coefficient of $f_n$ is
	$$
	a_k=\int_{[0,1)}f_n(x){\rm w}(k, x)\,dx= \int_{[0,1)} \sum_{m=0}^{2^n-1}b_m \xi_{m,n}(x) {\rm w}(k, x)\,dx
	$$
	$$
	=
	\sum_{m=0}^{2^n-1}b_m \int_{\Delta_{m,n}}{\rm w}(k, x)\,dx =
	\sum_{m=0}^{2^n-1}b_m \frac{1}{2^n}\omega^{n}_{k,m},
	$$   
	where $\omega^{n}_{k,m}$ is a value of ${\rm w}(k,\cdot)$ on $\Delta_{m,n}.$    
Let us denote $c_k :=b_k 2^{-n/2}$, $\bar{\omega}^{n}_{k,m}:=\omega^{n}_{k,m} 2^{-n/2}.$
Then 
$
a_k=\sum_{m=0}^{2^n-1}c_m \bar{\omega}^{n}_{k,m}, 
$
that is $a=c {\rm W}.$ 
Thus, $c$ is the fast Walsh-Fourier transform of $a.$

If $\tilde{x}_n$ minimizes the functional $\int_{\mathbb{R}_+}(x\oplus \tilde{x})^2 |f_n(x)|^2\,dx$ then 
$\tilde{x}_n$ can not be outside the support of $f_n.$ So, $\tilde{x} \in [0,\,1)=\cup_{k=\overline{0,2^n-1}}\Delta_{k,n}.$
Then, for $\tilde{x}\in \Delta_{k_0,n}$, we have 
$$
\int_{\mathbb{R}_+}(x\oplus \tilde{x})^2 |f_n(x)|^2\,dx=
\int_{[0,\,1)}(x\oplus \tilde{x})^2 \left|\sum_{k=0}^{2^n-1}b_k \xi_{k,n}(x)\right|^2\,dx
$$
$$
=\int_{[0,\,1)}(x\oplus \tilde{x})^2 \sum_{k=0}^{2^n-1}b_k^2 \xi_{k,n}(x)\,dx
=\sum_{k=0}^{2^n-1}b_k^2 \int_{\Delta_{k,n}}(x\oplus \tilde{x})^2 \,dx
$$
$$
=\left.\sum_{k=0}^{2^n-1}b_k^2\frac{x^3}{3}\right|_{\Delta_{k,n}\oplus \tilde{x}}
=\left.\sum_{k=0}^{2^n-1}b_{k\oplus k_0}^2 \frac{x^3}{3}\right|_{\Delta_{k,n}}=
\sum_{k=0}^{2^n-1}c_{k\oplus k_0}^2\frac{3k^2+3k+1}{3\times  2^{2n}}.
$$
  So, recalling Definition \ref{def1}, we get
$$
V(f_n):=\frac{1}{\|f_n\|^2_{L_2(\mathbb{R}_+)}} 
\min_{\tilde{x}}\int_{\mathbb{R}_+}(x\oplus \tilde{x})^2 |f(x)|^2\,dx=
\frac{1}{\sum_{k=0}^{2^n-1}|a_k|^2} 
\min_{k_0=\overline{0,2^n-1}}
\sum_{k=0}^{2^n-1}c_{k\oplus k_0}^2\frac{3k^2+3k+1}{3\times  2^{2n}}.
$$

The Walsh-Fourier transform of $f_n$ is
\begin{equation}
\widehat{f}_n(t)=
\sum_{k=0}^{2^n-1}a_k \int_{[0,\,1)}{\rm w}(x,\,t) {\rm w}(x,\,k) \, dx 
=
\sum_{k=0}^{2^n-1}a_k \chi_{[k,\,k+1)}(t).
\label{fWFfn}
\end{equation}
Then repeating the above calculations, we have
$$
V(\widehat{f}_n):=\frac{1}{\|\widehat{f}_n\|^2_{L_2(\mathbb{R}_+)}} \min_{\tilde{t}}\int_{\mathbb{R}_+}(t\oplus \tilde{t})^2 |\widehat{f}(t)|^2\,dt=
\frac{1}{\sum_{k=0}^{2^n-1}|c_k|^2} 
\min_{k_1=\overline{0,2^n-1}}
\sum_{k=0}^{2^n-1}a_{k\oplus k_1}^2\frac{3k^2+3k+1}{3}.
$$

To conclude the proof, it remains to show that 
$
UC_d(f)\,=\,\lim_{n\to \infty}UC_d(f_n).
$ 
We denote 
$V_0(g):=
\|g\|^2_{L_2(\mathbb{R}_+)} V(g)
 =\min_{\tilde{x}}\int_{\mathbb{R}_+}(x\oplus \tilde{x})^2 |g(x)|^2\,dx.$

Firstly, we prove  $\lim_{n\to\infty} V_0(f_n)=V_0(f).$ 
Assume that the minimum of the functional $V_0(f_n)$ is achieved at the point $\widetilde{x}_n^*$, the minimum of the functional  $V_0(f)$ is achieved at the point $\widetilde{x}^*.$
The functions $f_n$  converge uniformly on $[0,1)$ to $f$, i.e. for all $\varepsilon>0$ there exists $N\in{\mathbb N}$ such that for all $n\ge N$ and for all $  x\in[0,1)$ we have 
$\left||f(x)|-|f_n(x)|\right|\le|f(x)-f_n(x)|<\varepsilon.$ 
Then 
$$
|f(x)|^2-|f_n(x)|^2\le 2|f_n(x)||f(x)-f_n(x)|+|f(x)-f_n(x)|^2 \le
2|f_n(x)|\varepsilon + \varepsilon^2 \le
2(|f(x)|+\varepsilon)\varepsilon + \varepsilon^2.$$
After multiplication by $(x\oplus y)^2$ and integration over $[0,1)$ both sides of the above inequality, 
for all $y\in[0,1)$ and for all $n\ge N$ we get
$$
\int_{[0,1)}(x\oplus y)^2|f(x)|^2 dx-
\int_{[0,1)}(x\oplus y)^2|f_n(x)|^2 dx
\le  \varepsilon C
$$
where $C=\max\limits_{y\in[0,1)}\int_{[0,1)}
(x\oplus y)^2(2|f(x)|+
3\varepsilon)dx.$
The last inequality should be valid for $y=\widetilde{x}_n^*$ 
$$
\int_{[0,1)}(x\oplus \widetilde{x}_n^*)^2|f(x)|^2 dx-V_0(f_n)
\le  \varepsilon C \quad \forall n\ge N.
$$
Finally, we can decrease the left-hand side of the inequality by taking minimum of the functional over $\widetilde{x}_n^*$
$$
V_0(f)-V_0(f_n)
\le  \varepsilon C.
$$
Similarly, we can prove  the following inequality 
$$
V_0(f_n)-V_0(f)
\le  \varepsilon C.
$$
But it requires to start with 
$$
|f_n(x)|^2-|f(x)|^2\le 2|f(x)||f(x)-f_n(x)|+|f(x)-f_n(x)|^2 \le
2|f(x)|\varepsilon + \varepsilon^2 \quad \forall n\ge N, \forall x\in[0,1)
$$ 
and after the integration
take $y=\widetilde{x}^*.$ 
As a result, we get $\lim_{n\to\infty} V_0(f_n)=V_0(f).$

Now, let us prove $\lim_{n\to\infty} V_0(\widehat f_n)=V_0(\widehat f).$ 
Assume that the minimum of the functional $V_0( \widehat f_n)$ is achieved at the point $\widetilde{t}_n^*$, the minimum of the functional  $V_0(\widehat f)$ is achieved at the point $\widetilde{t}^*.$
By~(\ref{fWFfn}) we conclude that 
$|\widehat f_{n+1}(t)|^2\ge |\widehat f_{n}(t)|^2$ for all $t\in {\mathbb R}_+.$
After multiplication by $(t\oplus y)^2$ and integration over 
${\mathbb R}_+$ both sides of the above inequality, we get
$$\int_{{\mathbb R}_+}(t\oplus y)^2|\widehat f_{n+1}(t)|^2 dt\ge
\int_{{\mathbb R}_+}(t\oplus y)^2|\widehat f_{n}(t)|^2 dt\quad  \forall y\in {\mathbb R}_+.
$$
Thus, the last inequality should be valid for $y=\widetilde{t}_{n+1}^*$
$$V_0(\widehat f_{n+1})=
\int_{{\mathbb R}_+}(t\oplus \widetilde{t}_{n+1}^*)^2|\widehat f_{n+1}(t)|^2 dt\ge
\int_{{\mathbb R}_+}(t\oplus \widetilde{t}_{n+1}^*)^2|\widehat f_{n}(t)|^2 dt\ge 
V_0(\widehat f_{n}).
$$
Therefore, $V_0(\widehat f_{n+1})\ge V_0(\widehat f_{n})$ for all $n\in{\mathbb N},$ in particularly, $V_0(\widehat f)\ge V_0(\widehat f_{n})$.
Let us consider the difference
\begin{eqnarray*}
V_0(\widehat f)- V_0(\widehat f_{n}) & = & \min_{\widetilde t}
\int_{{\mathbb R}_+}(t\oplus \widetilde{t})^2|\widehat f(t)|^2 dt-
\int_{{\mathbb R}_+}(t\oplus \widetilde{t}_n^*)^2|\widehat f_n(t)|^2 dt
\\
& \le &
\int_{{\mathbb R}_+}(t\oplus \widetilde{t}_n^*)^2
\left(|\widehat f(t)|^2-|\widehat f_n(t)|^2\right)dt=
\int_{{\mathbb R}_+}(t\oplus \widetilde{t}_n^*)^2
\sum\limits_{k=2^n}^{\infty}|a_k|^2\chi_{[k,k+1)}(t)dt.
\end{eqnarray*}
 There exists $N\in{\mathbb N}$
such that $\widetilde{t}_n^*\in [0,2^N)$ and 
$\widetilde{t}^*\in [0,2^N)$ for all $n\in{\mathbb N}$ simultaneously. It can be shown by contradiction.
Indeed, assume that for any $N\in{\mathbb N}$ there exists $m>N$ such that $\widetilde{t}_m^*\ge 2^N.$
Then the following inequalities
\begin{eqnarray*}
V_0(\widehat f)\ge V_0(\widehat f_m) &=&
\int_{[0,2^N)}(t\oplus \widetilde{t}_m^*)^2|\widehat f_m(t)|^2\, dt+
\int_{[2^N,2^m)}(t\oplus \widetilde{t}_m^*)^2|\widehat f_m(t)|^2\, dt
\\
 & \ge & 
\int_{[0,2^N)}(t\oplus \widetilde{t}_m^*)^2|\widehat f_m(t)|^2\, dt
\ge 
2^N \sum_{k=0}^{2^N-1}|a_k|^2
\end{eqnarray*}
should be valid for all $N$. This leads to a contradiction.
The function 
$\int_{{\mathbb R}_+}(t\oplus y)^2
|\widehat f(t)|^2dt$ is bounded on $[0,2^N)$ (see Remark~\ref{minimum}).
Therefore, for all $\varepsilon>0$ there exists $M$ such that for all
$m>M,$ $m\in{\mathbb N}$ 
$$
\int_{[m,+\infty)}(t\oplus y)^2
|\widehat f(t)|^2dt=
\int_{{\mathbb R}_+}(t\oplus y)^2
\sum\limits_{k=m}^{\infty}|a_k|^2\chi_{[k,k+1)}(t)dt<\varepsilon.
$$
Then for $n$ such that $2^n>m$ we have
$V_0(\widehat f)- V_0(\widehat f_{n})<\varepsilon.$ 
Hence, $\lim_{n\to\infty} V_0(\widehat f_n)=V_0(\widehat f).$
Together with $\lim_{n\to\infty} V_0(f_n)=V_0(f),$ 
$\lim_{n \to \infty}\|f_n\|^2_{L_2(\mathbb{R}_+)}= \|f\|^2_{L_2(\mathbb{R}_+)}$, 
and 
$\lim_{n \to \infty}\|\widehat{f}_n\|^2_{L_2(\mathbb{R}_+)}= \|\widehat{f}\|^2_{L_2(\mathbb{R}_+)}$
we get the required statement for $UC_d.$
$\Box$

\begin{rem}
\label{rem4} 
It is easy to extend Lemma \ref{th1} to the functions of the form

\noindent
$g(x):=\chi_{[0,\,2^N)}(x) \sum_{k=0}^{\infty}a_k {\rm w}_k(x/2^N).$ Indeed, let
$g_n(x):=\chi_{[0,\,2^N)}(x) \sum_{k=0}^{2^n-1}a_k {\rm w}_k(x/2^N)$ be a partial sum of the above function $g,$  
 $f_n(x)=g_n(2^N x)$  the function defined in Lemma \ref{th1}. Then standard calculations show that
$\|g_n\|^2_2 = 2^N \|f_n\|^2_2,$ $\|\widehat{g}_n\|^2_2 = 2^N \|\widehat{f}_n\|^2_2,$
$\int_{\mathbb{R}_+}(x\oplus \tilde{x})^2 |g_n(x)|^2 \, dx = 
2^{3N}\int_{\mathbb{R}_+}(x\oplus (\tilde{x} 2^{-N}))^2 |f_n(x)|^2 \, dx$ and
$\int_{\mathbb{R}_+}(t\oplus \tilde{t})^2 |\widehat{g}_n(t)|^2 \, dt = 
2^{-N}\int_{\mathbb{R}_+}(t\oplus (\tilde{t} 2^{N}))^2 |\widehat{f}_n(t)|^2 \, dt.$
Hence, recalling the definition of $UC_d$ we get
$UC_d(g_n)=UC_d(f_n).$
The class of the functions of the form $g$ is rather large and important as  any orthogonal compactly supported dyadic scaling and wavelet 
functions belong to this set (see \cite[section 5]{FP}). 
\end{rem}

We denote $q_k:=\frac{3k^2+3k+1}{3\times  2^{n}}$ and suppose  $\|a\|=1,$ then
 $\|c\|=\|a W\|=1$ and the $UC_d(f_n)$ takes the form
\begin{equation*}
\label{UCnorm}
	UC_d(f_n)=\min_{k_1=\overline{0,2^n-1}}\sum_{k=0}^{2^n-1}a_{k\oplus k_1}^2 q_k \  \min_{k_0=\overline{0,2^n-1}}\sum_{k=0}^{2^n-1}c_{k\oplus k_0}^2 q_k.
\end{equation*}

Let us fix $n$. It follows from  (\ref{simp}) that the minimization problem 

\begin{equation}
\label{min1}
	\left\{
	\begin{array}{l}
	UC_d(f_n)\rightarrow \min \\
	\|a\|=1 
	\end{array}
	\right.
\end{equation}

 is equivalent to the following one

\begin{equation*}
\label{min2}
	\left\{
	\begin{array}{l}
	\sum_{k=0}^{2^n-1}a_{k}^2 q_k \  
	\sum_{k=0}^{2^n-1}c_{k}^2 q_k\rightarrow \min \\
	\|a\|=1 
	\end{array}
	\right.
\end{equation*}
 Using Wolfram Mathematica 8.0  
we solve numerically the last minimization problem for $n=2;\,3;\,4;\,5;\,6.$ The result is demonstrated in Table \ref{tab3}. 

\begin{table}[ht]
\centering
\caption{$UC_d(f_n)$}
\begin{tabular}{cccccc}
 \hline
$n$ & 2  & 3 & 4 & 5 & 6 \\
\hline
$\min_{f_n} UC_d(f_n)$ & 0.0891& 0.0882 & 0.0873 & 0.0881 & 0.0872\\
 \hline
\end{tabular}
\label{tab3}
\end{table}

\section{Examples}
\label{examsec}

\subsection{Lang's wavelet and scaling function}
To examine and illustrate the definition of the dyadic UC we use the first nontrivial example of orthogonal wavelets on the  Cantor dyadic group (see \cite {Lang}) . The dyadic scaling function is defined by
$$
\varphi_a(x)=\frac12 \chi_{\left[0,\,1\right)}\left(\frac x2\right) \left(1+a\sum_{j=0}^{\infty}b^j {\rm w}\left(2^{j+1}-1,\, \frac x2\right)\right),
\quad
 \widehat{\varphi}_a=\chi_{\left[0,\,1/2\right)}+a \sum_{j=0}^{\infty}b^j \chi_{[2^j-1/2,\,2^j)},
$$
 where  $0<a\leq 1,$  $a^2+b^2=1,$ $a,b\in {\mathbb R}.$
 The corresponding wavelet is defined by
$$
\psi_a(x)\,= \, 2 a_0 \varphi_a(2x\oplus 1)-2 a_1 \varphi_a(2x)
+2 a_2 \varphi_a(2x\oplus 3)-2 a_3 \varphi_a(2x\oplus 2),
$$
where
$
a_0\,=\,(1+a+b)/4,\quad a_1\,=\,(1+a-b)/4,\quad a_2\,=\,(1-a-b)/4,\quad a_3\,=\,(1-a+b)/4.
$
Then the wavelet system $\{\psi_{j,k}\}_{j\in\mathbb Z, k\in \mathbb R_+}$ forms an orthonormal basis in $L_2(\mathbb{R}_+).$

The integrals defying the dyadic UC for the scaling and wavelet functions  are
	\begin{eqnarray*}
\int_{\mathbb{R}_+}(x \oplus \tilde{x})^2|\varphi_a(x)|^2\,d x
 = 
	\frac {4}3 +\frac 14 \mathrm{w}\left(1,\,\frac {\tilde{x}}{2}\right)\left(-4a+ a\, b\, \mathrm{w}(1,\,\tilde{x})\right)\phantom{123456789012345}\\
	-
	\frac {a^2 b}2 \sum\limits_{j=0}^{\infty} 
	\left(\frac {b^2}{2}\right)^j \mathrm{w}(2^j, \, \tilde{x})
	+
	\frac {a^2 b^2}{16} \sum\limits_{j=0}^{\infty} 
	\left(\frac {b^2}{4}\right)^j \mathrm{w} (2^j\oplus2^{j+1},\,\tilde{x});\\
	\int_{\mathbb{R_+}}(t \oplus \tilde{t})^2|\widehat{\varphi_a}(t)|^2\,d t
	= A(0,\,\tilde{t})
	+
	a^2 \sum\limits_{j=0}^{\infty} 
	b^{2j} 	A(2^j-1/2,\,\tilde{t}),\phantom{1234567890123456}
	\end{eqnarray*}

\begin{eqnarray*}
\int_{\mathbb{R}_+}(x \oplus \tilde{x})^2|\psi_a(x)|^2\,d x
 = 
\frac 43 -a \mathrm{w}_1\left(\frac {\tilde{x}}2\right) 
	- \frac {ab}4  \mathrm{w}_1\left(\frac {\tilde{x}}2\right) \mathrm{w}(1,\, \tilde{x})
	\phantom{12345678901234567}
	\\
	-\frac {a^2 b}2\left(-\mathrm{w}(1,\,\tilde{x}) + \frac b8 \mathrm{w}(3,\,\tilde{x})\right)
	+
	a^2\left(\frac 14 +\frac {a^2}4\right) 
	\left(-b
	\sum\limits_{j=0}^{\infty} 
	\left(\frac {b^2}2\right)^j \mathrm{w}(2^{j+1},\,\tilde{x})\right.
	\\   
	\left.
	+\frac {b^2}{16} \sum\limits_{j=0}^{\infty} 
	\left(\frac {b^2}4\right)^j \mathrm{w} (2^{j+1}\oplus2^{j+2},\,\tilde{x})\right)
	+\frac {a^3}4 \mathrm{w}\left(1,\,\frac {\tilde{x}}2\right) b \sum\limits_{j=0}^{\infty} 
	\left(\frac {b^2}2\right)^j \mathrm{w} (2^{j+1},\,\tilde{x});
	\\
	\int_{\mathbb{R}_+}(t \oplus \tilde{t})^2|\widehat{\psi}_a(t)|^2\,d t
	\phantom{12345678901234567890123456789012345678901234567890}\\
	=
	b^2
A(1/2,\,\tilde{t})
	+
	a^2 \sum\limits_{j=0}^{\infty} 
	b^{2j} 	A(2^j-1,\,\tilde{t})+
	a^4 \sum\limits_{j=1}^{\infty} 
	b^{2j} 	A(2^j-1/2,\,\tilde{t}),
	\end{eqnarray*}
	where $A(\xi,\,\eta)=\frac 13 ((\inf \{[\xi,\,\xi+1/2)\oplus \eta\} )+ 1/2)^3-\frac 13 (\inf \{[\xi,\,\xi+1/2)\oplus \eta\} )^3.$	
It turns out that
	\begin{gather*}
  UC_{d}(\varphi_a),\,  UC_{d}(\psi_a)<\infty  \Leftrightarrow \sqrt{3}/2 <a\leq 1.
\end{gather*} 
	The dyadic UCs for the different values of the parameter $a$ are collected in  Table \ref{tab2} and Table \ref{tabpsi}.
	The best localized function here is the Haar scaling function. It corresponds to the case $a=1$. 
	
\begin{table}[ht]
\centering
\caption{The dyadic uncertainty constants for  $\varphi_a$.}
\begin{tabular}{cccccc} 
\hline
$a$ & $V(\varphi_a)$ & $\tilde{x}_{0}(\varphi_a)$ & $V(\widehat{\varphi}_a)$ & $\tilde{t}_{0}(\varphi_a)$ & $UC_d(\varphi_a)$\\ 
\hline
$0.9$  &  $0.346$ & $0$ & $1.29$ & $[1/2,\,1)$ & $0.446$ 
 \\ 
$0.95$ &  $0.315$ & $0$ & $0.482$ & $[1/2,\,1)$ & $0.152$
\\
1 &  $1/3$ & $[0,\,1)$ & $1/3$ & $[0,\,1)$ &  $1/9$
\\
 \hline
\end{tabular}
\label{tab2}
\end{table}

\begin{table}[ht]
\centering
\caption{The dyadic uncertainty constants for  $\psi_a$.}
\begin{tabular}{cccccc} 
\hline
$a$ & $\Delta_d^2(\psi_a)$ & $\tilde{x}_{0}(\psi_a)$ & $\Delta_d^2(\widehat{\psi}_a)$ & $\tilde{t}_{0}(\varphi_a)$ & $UC_d^2(\psi_a)$ \\ 
\hline
$0.9$  &  $0.280$ & $0.5$ & $7.438$ & $[3/2,\,2)$ & $2.083$ 
 \\ 
$0.95$ &  $0.254$ & $ 0.5$ & $1.546$ & $[3/2,\,2)$ & $0.393$
\\
1 &  $1/3$ & $[0,\,1)$ & $1/3$ & $[0,\,1)$ &  $1/9$
\\
 \hline
\end{tabular}
\label{tabpsi}
\end{table}

\subsection{Dyadic wavelet frames with good localization}
\begin{enumerate}
	\item 

Let us consider generators of normalized tight frames~\cite[Example 3.2]{F12} for $L_2({\mathbb R}_+):$
$$g_{l,s}(x)=2^{-s}\chi_{[0,2^s)}\mathrm{w}(l,\,2^{-s}x),$$
where $l\in {\mathbb N},$ $s\in {\mathbb Z}_+.$
The Walsh-Fourier transform of $g_{l,s}$ is $\widehat{g_{l,s}}=\chi_{U_{l,s}},$
where $U_{l,s}=2^{-s}(l\oplus[0,1)).$
Suppose that $\psi=g_{l,s}.$ Then $\{\psi_{j,\alpha}\}$ is a normalized tight frame for $L_2({\mathbb R}_+)$. 
For all $l\in {\mathbb N}, s\in {\mathbb Z}_+$ the dyadic UC is $UC_d(g_{l,s})=\frac 19.$ 

%

\item

As it was noted in Table \ref{tab3}
numerically $\min UC_d(f_n)\simeq 0.0891$ for $n=2.$
Let us try to find a frame generator such that its dyadic UC is close to this
 value. Let $\psi=\chi_{[0,\,1)}(x) \sum_{k=0}^{3} a_k {\rm w}(k,\,x).$ From the frame criteria, we should provide zero moment for the frame
  generator $\psi$ or, equivalently, $\widehat\psi(0)=0.$ Thus, we assume that $a_0=0.$  Using Wolfram Mathematica 8.0  
we solve numerically the  minimization problem
(\ref{min1}).
The coefficients are
$(a_0,a_1,a_2,a_3)=(0, 0.094206, 0.551564, 0.828796).$ 
Using Theorem 3.2 in \cite{F12}, we compute the frame bounds for the frame $\{\psi_{jk}\},$ namely $A=0.313098,$ $B=0.695777.$
The dyadic UC is $UC_d(\psi)=0.091286$ and it is close to the minimal possible constant for $n=2.$

The same computations can be done for the case $n=3.$
Let $\psi=\chi_{[0,\,1)}(x) \sum_{k=0}^{7} a_k {\rm w}(k,\,x).$ The  minimum for  $UC_d(\psi)$ is delivered by the coefficients 
$(a_0,a_1,a_2,a_3,a_4,a_5,a_6,a_7)=
(0, 0.001335, -0.009155, -0.022170, -0.067567, 
-0.138436, -0.601657, -0.783391).$ The frame bounds for the frame $\{\psi_{jk}\}$ are $A=0.004649,$ $B=0.614194.$ The dyadic UC is $UC_d(\psi)=0.0882147.$

\end{enumerate}

\section*{Acknowledgments}
The work  is supported by the RFBR-grant \#12-01-00216 and 
  by the grant of President RF \#MK-1847.2012.1.

\end{document}